\documentclass[10pt]{article}
\usepackage[centertags]{amsmath}
\usepackage{amsfonts}
\usepackage{amssymb}
\usepackage{amsthm}
\usepackage{newlfont}
\usepackage{graphicx}

\newfont{\bb}{msbm10 at 12pt}

\input{tcilatex}
\begin{document}

\title{\textbf{On the natural lift curves for the Involute spherical
indicatrices in Minkowski 3-space}}
\author{Mustafa Bilici\thanks{%
Corresponding author. Tel.: +903623121919, E-mail adress: } \\
Ondokuz May\i s University, Education Faculty, Department of Mathematics,\\
55200 Kurupelit, Samsun, Turkey \ \ \ \ \  \and Ahmad T. Ali \\
King Abdul Aziz University, Faculty of Science, Department of Mathematics, \\
PO Box 80203, Jeddah, 21589, Saudi Arabia.}
\maketitle

\begin{abstract}
The aim of this paper is to determine criteria of being integral curve for
the geodesic spray of the natural lift curves of the spherical indicatrices
of the involutes of a given spacelike curve $\alpha $ with a timelike
binormal in Minkowski 3-space $\hbox{\bf E}_{1}^{3}$. Furthermore, some
interesting results about the spacelike evolute curve with timelike binormal
and spacelike or timelike Darboux vector $\omega ~$were obtained, depending
on the assumption that the natural lift curves $\overline{\alpha ^{\ast }}%
_{t^{\ast }},~\overline{\alpha ^{\ast }}_{n^{\ast }}~$and $\overline{\alpha
^{\ast }}_{b^{\ast }}$ of the spherical indicatrices $\alpha _{t^{\ast
}}^{\ast }$ ,$\alpha _{n^{\ast }}^{\ast }$ and $\alpha _{b^{\ast }}^{\ast }$
of the involute curve $\alpha ^{\ast }$ should be the integral curve on the
tangent bundle $T\left( S_{1}^{2}\right) $ or $T\left( H_{0}^{2}\right) $ .
Additionally we illustrate an example of our main results.

\textbf{M.S.C. 2000}: 53C40, 53C50.\newline
\end{abstract}

\section{Introduction}

One of the most significant curve is an involute of a given curve. C.
Huygens, who is also known for his works in optics, discovered involutes
while trying to build a more accurate clock. The original curve is called an
evolute. A curve can have any number of involutes, thus a curve is an
evolute of each of its involutes and an involute of its evolute. The normal
to a curve is tangent to its evolute and the tangent to a curve is normal to
its involutes. In addition to this, involute-evolute curve couple is a well
known concept in the classical differential geometry, see \cite{mill, hs,
cal}. The basic local theory of space curve are mainly developed by the
Frenet-Serret theorem which expresses the derivative of a geometrically
chosen basis of \textit{$\hbox{\bf E}^{3}$} by the aid of itself is proved.
Then it is observed that by the solution of some of special ordinary
differential equations, further classical topics, for instance spherical
curves, Bertrand curves, involutes and evolutes are investigated, see for
details \cite{Do Carmo}.

In differential geometry, especially the theory of space curves, the Darboux
vector is the areal velocity vector of the Frenet frame of a space curve. It
is named after Gaston Darboux who discovered it. In terms of the
Frenet-Serret apparatus, the Darboux vector $\omega ~$can be expressed as $%
\omega =\tau t+\kappa b$. In addition to this, the concepts of the natural
lift and the geodesic sprays have been given by Thorpe in 1979 \cite{thorpe}%
. \c{C}al\i \c{s}kan et al. \cite{call} have studied the natural lift curves
and the geodesic sprays in the Euclidean 3-space \textit{$\hbox{\bf E}^{3}$}%
. Then Bilici et al. \cite{bilici1} have proposed the natural lift curves
and the geodesic sprays for the spherical indicatrices of the the
involute-evolute curve couple in \textit{$\hbox{\bf E}^{3}$}.

Spherical images (indicatrices) are a well known concept in classical
differential geometry of curves \cite{Do Carmo}. Kula and Yayl\i\ \cite{Kula}
have studied spherical images of the tangent indicatrix and binormal
indicatrix of a slant helix and they have shown that the spherical images
are spherical helices. In recent years some of the classical differential
geometry topics have been extended to Lorentzian geometry. In \cite{Suha} S%
\"{u}ha at all investigated tangent and trinormal spherical images of
timelike curve lying on the pseudo hyperbolic space in Minkowski space-time. 
\.{I}yig\"{u}n \cite{Esen} defined the tangent spherical image of a unit
speed timelike curve lying on the on the pseudo hyperbolic space in $%
H_{0}^{2}$. In \cite{bilici2} author adapted this problem for the spherical
indicatrices of the involutes of a timelike curve in Minkowski 3-space $%
\hbox{\bf E}_{1}^{3}$. However, this problem is not solved in other cases of
the space curve.

In the present paper, the natural lift curves for the spherical indicatrices
of the involutes of a given spacelike curve with a timelike binormal have
been investigated in Minkowski 3-space $\hbox{\bf E}_{1}^{3}$. With this aim
we translate tangents of the involutes of a spacelike curve with a timelike
binormal curve to the center of the unit hypersphere $S_{1}^{2}$ we obtain a
spacelike curve $\alpha _{t^{\ast }}^{\ast }=t^{\ast }$ on the unit
hypersphere . This curve is called the first spherical indicatrix or tangent
indicatrix of $\alpha ^{\ast }$. One consider the principal normal
indicatrix $\overline{\alpha ^{\ast }}_{n^{\ast }}=n^{\ast }~$and the
binormal indicatrix $\overline{\alpha ^{\ast }}_{b^{\ast }}=b^{\ast }~$on
the unit hypersphere $H_{0}^{2}$. Then the natural lift curves of the
spherical indicatrices of the involutes of a given spacelike curve $\alpha $
with a timelike binormal are investigated in Minkowski 3-space $\hbox{\bf E}%
_{1}^{3}$ and some new results were obtained. We hope these results will be
helpful to mathematicians who are specialized on mathematical modeling.

\section{Preliminaries}

Let $M$ be a hypersurface in $\hbox{\bf E}_{1}^{3}$ equipped with a metric $%
g $, where the metric $g$ means a symmetric non-degenerate $\left(
0,2\right) $ tens\"{o}r field on $M$ with constant signature. For a
hypersurface $M$, let $TM$ be the set $\cup \left\{ T_{p}\left( M\right)
:~p\in M\right\} ~$of all tangent vectors to $M$. A technicality: For each $%
p\in M~$replace $0\in T_{p}\left( M\right) ~$by $0_{p}$ (other-wise the zero
tangent vector is in every tangent space). Then each $v\in TM~$is in a
unique $T_{p}\left( M\right) $, and the projection $\pi :TM\rightarrow M$
sends $v$ to $p$. Thus $\pi ^{-1}\left( p\right) =T_{p}\left( M\right) $.
There is a natural way to make $TM$ a manifold, called the \textit{tangent
bundle} of $M$.

A vector field $X\in \chi \left( M\right) ~$is exactly a smooth section of $%
TM$, that is, a smooth function $X:M\rightarrow TM~$such that $\pi \circ
X=I~\left( identity\right) $. Let $M$ be a hypersurface in $\hbox{\bf E}%
_{1}^{3}$. A curve $\alpha :I\rightarrow TM~$\ is an integral curve of $X\in
\chi \left( M\right) ~$provided $\alpha ^{\prime }=X_{\alpha }~$; that is,

\begin{equation*}
\frac{d}{ds}\left( \alpha \left( s\right) \right) =X\left( \alpha \left(
s\right) \right) ~\text{for all }s\in I,\text{~\cite{Onei},}
\end{equation*}%
For any parametrized curve $\alpha :I\rightarrow TM$, the parametrized curve
given by $\overline{\alpha }:I\rightarrow TM$

\begin{equation*}
s\rightarrow \overline{\alpha }\left( s\right) =\left( \alpha \left(
s\right) ,\alpha ^{\prime }\left( s\right) \right) =\alpha ^{\prime }\left(
s\right) \mid _{\alpha \left( s\right) }
\end{equation*}

is called the \textit{natural lift} of $\alpha ~$on $TM$. Thus, we can write

\begin{equation}
\frac{d\overline{\alpha }}{ds}=\frac{d}{ds}\left( \alpha ^{\prime }\left(
s\right) \mid _{\alpha \left( s\right) }\right) =D_{\alpha ^{\prime }\left(
s\right) }\alpha ^{\prime }\left( s\right) ,  \label{u1}
\end{equation}%
where $D$ is the standard connection on $\hbox{\bf E}_{1}^{3}$.

For $v\in TM~$, the smooth vector field $X\in \chi \left( M\right) ~$defined
by

\begin{equation}
X\left( v\right) =\varepsilon g\left( v,S\left( v\right) \right) \xi \mid
_{\alpha \left( s\right) },~\varepsilon =g\left( \xi ,\xi \right)  \label{u2}
\end{equation}%
is called the \textit{geodesic spray} on the manifold $TM$, where $\xi ~$is
the unit normal vector field of $M$ and $S$ is the shape operator of $M$.

The Minkowski three-dimensional space $\hbox{\bf E}_{1}^{3}$ is the real
vector space $\mathbb{R}^{3}$ endowed with the standard flat Lorentzian
metric given by \cite{ali3} 
\begin{equation*}
g=-dx_{1}^{2}+dx_{2}^{2}+dx_{3}^{2},
\end{equation*}%
where $(x_{1},x_{2},x_{3})$ is a rectangular coordinate system of $%
\hbox{\bf
E}_{1}^{3}$. If $u=(u_{1},u_{2},u_{3})$ and $v=(v_{1},v_{2},v_{3})$ are
arbitrary vectors in $\hbox{\bf E}_{1}^{3}$ then we define the Lorentzain
vector product of $u~$and $v~$as the following:

\begin{equation*}
u\times
v=(u_{3}v_{2}-u_{2}v_{3},u_{3}v_{1}-u_{1}v_{3},u_{1}v_{2}-u_{2}v_{1}).
\end{equation*}

Since $g$ is an indefinite metric, recall that a vector $v\in \hbox{\bf E}%
_{1}^{3}$ can have one of three Lorentzian characters: it can be space-like
if $g(v,v)>0$ or $v=0$, timelike if $g(v,v)<0$ and null if $g(v,v)=0$ and $%
v\neq 0$. Similarly, an arbitrary curve $\alpha =\alpha (s)$ in $\hbox{\bf E}%
_{1}^{3}$ can locally be spacelike, timelike or null (lightlike), if all of
its velocity vectors $\alpha ^{\prime }$ are respectively spacelike,
timelike or null (lightlike), for every $s\in I\subset \mathbb{R}$. The
pseudo-norm of an arbitrary vector $a\in \hbox{\bf E}_{1}^{3}$ is given by $%
\left\Vert a\right\Vert =\sqrt{\left\vert g(a,a)\right\vert }$. $\alpha $ is
called an unit speed curve if velocity vector $\sigma $ of $\alpha $
satisfies $\left\Vert \sigma \right\Vert =1$. For vectors $v,w\in 
\hbox{\bf
E}_{1}^{3}$ it is said to be orthogonal if and only if $g(v,w)=0$.\newline

Denote by $\left\{ \hbox{\bf t},\hbox{\bf n},\hbox{\bf b}\right\} $ the
moving Frenet frame along the curve $\alpha $ in the space $\hbox{\bf E}%
_{1}^{3}$. For an arbitrary curve $\alpha $ with first and second curvature, 
$\kappa $ and $\tau $ in the space $\hbox{\bf E}_{1}^{3}$, the following
Frenet formulae are given in \cite{ilarslan}: If $\alpha $ is a spacelike
curve with a timelike binormal vector $\hbox{\bf b}$, then the Frenet
formulae read

\begin{equation}
\left[ 
\begin{array}{c}
\hbox{\bf t}^{\prime } \\ 
\hbox{\bf n}^{\prime } \\ 
\hbox{\bf b}^{\prime }%
\end{array}%
\right] =\left[ 
\begin{array}{ccc}
0 & \kappa & 0 \\ 
-\kappa & 0 & \tau \\ 
0 & \tau & 0%
\end{array}%
\right] \left[ 
\begin{array}{c}
\hbox{\bf t} \\ 
\hbox{\bf n} \\ 
\hbox{\bf b}%
\end{array}%
\right] ,  \label{u3}
\end{equation}%
where $g(\hbox{\bf t},\hbox{\bf t})=1,\,\,g(\hbox{\bf n},\hbox{\bf n}%
)=1,\,\,g(\hbox{\bf b},\hbox{\bf b})=-1,\,\,g(\hbox{\bf t},\hbox{\bf n})=g(%
\hbox{\bf t},\hbox{\bf b})=g(\hbox{\bf n},\hbox{\bf b})=0.$

The angle between two vectors in Minkowski space is defined by \cite{rat}:

\textbf{Definition 2.1.} Let $X$ and $Y$ be spacelike vectors in $%
\hbox{\bf
E}_{1}^{3}$ that span a spacelike vector subspace, then we have $%
|g(X,Y)|\leq \Vert X\Vert \Vert Y\Vert $ and hence, there is a unique
positive real number $\theta $ such that 
\begin{equation*}
|g(X,Y)|=\Vert X\Vert \Vert Y\Vert \mathrm{cos}\theta .
\end{equation*}%
The real number $\theta $ is called the Lorentzian spacelike angle between $%
X $ and $Y$.

\textbf{Definition 2.2.} Let $X$ and $Y$ be spacelike vectors in $%
\hbox{\bf
E}_{1}^{3}$ that span a timelike vector subspace, then we have $%
|g(X,Y)|>\Vert X\Vert \Vert Y\Vert $ and hence, there is a unique positive
real number $\theta $ such that 
\begin{equation*}
|g(X,Y)|=\Vert X\Vert \Vert Y\Vert \mathrm{cosh}\theta .
\end{equation*}%
The real number $\theta $ is called the Lorentzian timelike angle between $X$
and $Y$.

\textbf{Definition 2.3.} Let $X$ be a spacelike vector and $Y$ a positive
timelike vector in $\hbox{\bf E}_{1}^{3}$, then there is a unique
non-negative real number $\theta $ such that 
\begin{equation*}
|g(X,Y)|=\Vert X\Vert \Vert Y\Vert \mathrm{sinh}\theta .
\end{equation*}%
The real number $\theta $ is called the Lorentzian timelike angle between $X$
and $Y$.

\textbf{Definition 2.4.} Let $X$ and $Y$ be positive (negative) timelike
vectors in $\hbox{\bf E}_{1}^{3}$, then there is a unique non-negative real
number $\theta $ such that 
\begin{equation*}
g(X,Y)=\Vert X\Vert \Vert Y\Vert \mathrm{cosh}\theta .
\end{equation*}%
The real number $\theta $ is called the Lorentzian timelike angle between $X$
and $Y$.

The Darboux vector for the spacelike curve with a timelike binormal is
defined by \cite{ugur}:

\begin{equation*}
\omega =\tau t-\kappa b.
\end{equation*}%
There are two cases corresponding to the causal characteristic of Darboux
vector\ $\omega $

\textbf{Case 1}\textit{. }If $\left\vert \kappa \right\vert <\left\vert \tau
\right\vert $, then $\omega ~$is a spacelike vector. In this situation, we
can write

\begin{equation*}
\kappa =\left\Vert \omega \right\Vert \sin h\theta ,~\tau =\left\Vert \omega
\right\Vert \cos h\theta ,~g\left( \omega ,\omega \right) =\left\Vert \omega
\right\Vert ^{2}=\tau ^{2}-\kappa ^{2}
\end{equation*}%
and the unit vector $c$ of direction $\omega $ is

\begin{equation*}
c=\frac{1}{\left\Vert \omega \right\Vert }\omega =\cos h\theta t-\sin
h\theta b,
\end{equation*}%
where $\theta ~$is the Lorentzian timelike angle between $-b$ and timelike
unit vector $c^{\prime }$ Lorentz orthogonal to the normalisation of the
Darboux vector $c$ as Fig. 1.

\ \ \ \ \ \ \ \ \ \ \ \ \ \ \ \ \ \ \ \ \ \ \ \ \ \ \ \ \ \ \ \ \FRAME{dtbpFU%
}{3.7343in}{3.6608in}{0pt}{\Qcb{\textbf{Figure 1. }Lorentzian timelike angle 
$\protect\theta $}}{}{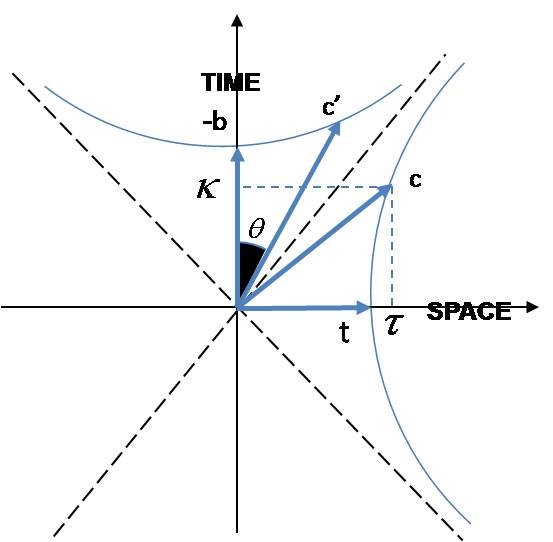}{\special{language "Scientific Word";type
"GRAPHIC";maintain-aspect-ratio TRUE;display "USEDEF";valid_file "F";width
3.7343in;height 3.6608in;depth 0pt;original-width 3.6867in;original-height
3.614in;cropleft "0";croptop "1";cropright "1";cropbottom "0";filename
'fig1.jpg';file-properties "XNPEU";}}

\textbf{Case 2}\textit{. }If $\left\vert \kappa \right\vert >\left\vert \tau
\right\vert $, then $\omega ~$is a timelike vector. In this situation, we
have

\begin{equation*}
\kappa =\left\Vert \omega \right\Vert \cos h\theta ,~\tau =\left\Vert \omega
\right\Vert \sin h\theta ,~g\left( \omega ,\omega \right) =-\left\Vert
\omega \right\Vert ^{2}=\kappa ^{2}-\tau ^{2}
\end{equation*}%
and the unit vector $c$ of direction $\omega $ is

\begin{equation*}
c=\frac{1}{\left\Vert \omega \right\Vert }\omega =\sin h\theta t-\cos
h\theta b,
\end{equation*}
\ 

\textbf{Proposition 2.5. }Let $\alpha ~$be a spacelike (or timelike) curve
with curvatures $\kappa $ and $\tau $. The curve is a general helix if and
only if $\frac{\tau }{\kappa }$ is constant, \cite{barr}.

\textbf{Remark 2.6 . }We can easily see from Lemma 3.2, 3.3, and 3.4 in \cite%
{ali1} that: $\frac{\tau (s)}{\kappa (s)}=\cot \theta $, $\frac{\tau (s)}{%
\kappa (s)}=\coth \theta $ or $\frac{\tau (s)}{\kappa (s)}=\tanh \theta $,
if $\theta =$constant then $\alpha $ is a general helix.

\textbf{Lemma 2.7. }\textit{The natural lift }$\overline{\alpha }$\textit{\
of the curve }$\alpha $\textit{\ is an integral curve of the geodesic spray }%
$X$\textit{\ if and only if }$\alpha $\textit{\ is a geodesic on }$M$\textit{%
\ }\cite{bilici2}.

\textbf{Remark 2.8. }Let $\alpha $ be a spacelike curve with a timelike
binormal. In this situation its involute curve $\alpha ^{\ast }$ must be a
spacelike curve with a spacelike or timelike binormal. $\left( \alpha
,~\alpha ^{\ast }\right) $ being the involute-evolute curve couple, the
following lemma was given by \cite{bili}.

\textbf{Lemma 2.9.}\textit{\ Let }$\left( \alpha ,~\alpha ^{\ast }\right) $%
\textit{\ be the involute-evolute curve couple. The relations between the
Frenet vectors of the curve couple as follow.}

\textbf{I}\textit{. }If $\omega ~$is a spacelike vector $\left( \left\vert
\kappa \right\vert <\left\vert \tau \right\vert \right) $, then

\begin{center}
$\left[ 
\begin{array}{c}
\hbox{\bf t}^{\ast } \\ 
\hbox{\bf n}^{\ast } \\ 
\hbox{\bf b}^{\ast }%
\end{array}%
\right] =\left[ 
\begin{array}{ccc}
0 & 1 & 0 \\ 
\sinh \theta & 0 & -\cosh \theta \\ 
-\cosh \theta & \tau & \sinh \theta%
\end{array}%
\right] \left[ 
\begin{array}{c}
\hbox{\bf t} \\ 
\hbox{\bf n} \\ 
\hbox{\bf b}%
\end{array}%
\right] .$
\end{center}

\textbf{II}\textit{. }If $\omega ~$is a timelike vector $\left( \left\vert
\kappa \right\vert >\left\vert \tau \right\vert \right) $, then

\begin{center}
$\left[ 
\begin{array}{c}
\hbox{\bf t}^{\ast } \\ 
\hbox{\bf n}^{\ast } \\ 
\hbox{\bf b}^{\ast }%
\end{array}%
\right] =\left[ 
\begin{array}{ccc}
0 & 1 & 0 \\ 
-\cosh \theta & 0 & \sinh \theta \\ 
-\sinh \theta & \tau & \cosh \theta%
\end{array}%
\right] \left[ 
\begin{array}{c}
\hbox{\bf t} \\ 
\hbox{\bf n} \\ 
\hbox{\bf b}%
\end{array}%
\right] .$
\end{center}

\textbf{Remark 2.10.}\textit{\ }In this situation I.,\ the causal
characteristics of the Frenet frame of the involute curve $\alpha ^{\ast }$
is $\left\{ t^{\ast }~spacelike,~n^{\ast }~timelike,~b^{\ast
}~spacelike\right\} $.\textit{\ }If $\alpha $ is a spacelike curve with
timelike $\omega $, then the causal characteristics of the Frenet frame of
the curve $\alpha ^{\ast }$ must be of he form $\left\{ t^{\ast
}~spacelike,~n^{\ast }~spacelike,~b^{\ast }~timelike\right\} $.\textit{\ }

\textbf{Definition 2.10}. Let $S_{1}^{2}$ and $H_{0}^{2}$ be hyperspheres in 
$\hbox{\bf E}_{1}^{3}$. The Lorentzian sphere and hyperbolic sphere of
radius 1 in are given by

\begin{equation*}
S_{1}^{2}=\left\{ a=(a_{1},a_{2},a_{3})\epsilon \hbox{\bf E}_{1}^{3}:g\left(
a,a\right) =1\right\}
\end{equation*}%
\qquad

and

\qquad 
\begin{equation*}
H_{0}^{2}=\left\{ a=(a_{1},a_{2},a_{3})\epsilon \hbox{\bf E}_{1}^{3}:g\left(
a,a\right) =-1\right\}
\end{equation*}

respectively, \cite{Onei}.

\section{\textbf{The natural lift curves for the spherical indicatrices of
the involutes of a spacelike curve with a timelike binormal}}

\subsection{\textbf{The natural lift} of tangent indicatrix \textbf{of the} 
\textbf{curve }$\protect\alpha ^{\ast }$}

Let $\alpha ~$be a spacelike curve with timelike binormal and spacelike $%
\omega $ $\left( \left\vert \kappa \right\vert <\left\vert \tau \right\vert
\right) $. We will investigate how evolute curve $\alpha $ must be a curve
satisfying the condition that the natural lift curve $\overline{\alpha
^{\ast }}_{t^{\ast }}$ is an integral curve of geodesic spray, where $\alpha
_{t^{\ast }}^{\ast }$ is the spherical indicatrix of tangent vector of
involute curve $\alpha ^{\ast }$.

If the natural lift curve $\overline{\alpha ^{\ast }}_{t^{\ast }}$\ is an
integral curve of the geodesic spray, then by means of Lemma 2.1.\qquad
\qquad 
\begin{equation}
\overline{D}_{\alpha _{t^{\ast }}^{\prime \ast }}\alpha _{t^{\ast }}^{\prime
\ast }=0,  \label{u4}
\end{equation}%
where $\overline{D}$ is the connection on the Lorentzian sphere $S_{1}^{2}$
and the equation of the spherical indicatrix of tangent vector of the
involute curve $\alpha ^{\ast }~$is $\alpha _{t^{\ast }}^{\ast }=t^{\ast }$.
Thus from Lemma 2.2.I and the last equation we obtain

\begin{equation*}
-\frac{\theta ^{\prime }}{\left\Vert \omega \right\Vert }\cosh \theta +\frac{%
\theta ^{\prime }}{\left\Vert \omega \right\Vert }\sinh \theta =0.
\end{equation*}%
Because of $\left\{ t,~n,~b\right\} $ are linear independent, we can easily
see that 
\begin{equation*}
\theta =\text{cons}\tan \text{t},
\end{equation*}%
according to Remark 2.1, we have

\begin{equation*}
\frac{\tau }{\kappa }=\coth \theta =\text{cons}\tan \text{t}.
\end{equation*}

\textbf{Result 3.1.1. }If the curve $\alpha $ is a general helix, then the
spherical indicatrix $\alpha _{t^{\ast }}^{\ast }$ of the involute curve $%
\alpha ^{\ast }$ is a geodesic on the Lorentzian sphere $S_{1}^{2}$. In this
case, from the Lemma 2.1 the natural lift $\overline{\alpha ^{\ast }}%
_{t^{\ast }}$ of $\alpha _{t^{\ast }}^{\ast }$ is an integral curve of the
geodesic spray on the tangent bundle $T\left( S_{1}^{2}\right) $. In the
case of a spacelike curve with timelike binormal and timelike $\omega $,
similar result can be easily obtained in following same procedure.\ \ \ \ \
\ \ \ \ \ \ \ \ \ \ \ \ \ \ \ \ 

\textbf{Remak 3.1.2. }From the classification of all W-curves (i.e. a curves
for which a curvature and a torsion are constants) in \cite{ali1, walrave},
Case 1. and\ Case 2.\ we have following results with relation to curve $%
\alpha $. \ \ \ \ 

\textbf{Result 3.1.3. }If the curve $\alpha $ with $\kappa =$cons$\tan $t $%
>0,~\tau =$cons$\tan $t$\neq 0$ and $\kappa <\left\vert \tau \right\vert $
then $\alpha $ is a part of a spacelike hyperbolic helix,

\begin{equation*}
\alpha \left( s\right) =\frac{1}{\left\Vert \omega \right\Vert ^{2}}\left(
\kappa \sinh \left[ \left\Vert \omega \right\Vert s\right] ,~\kappa \cosh %
\left[ \left\Vert \omega \right\Vert s\right] ,~\tau \left\Vert \omega
\right\Vert s\right) .
\end{equation*}

\textbf{Result 3.1.4. }Let $\alpha $ be a spacelike curve with timelike
binormal and timelike $\omega $.$~$If the curve $\alpha $ with $\kappa =$cons%
$\tan $t $>0,~\tau =$cons$\tan $t$\neq 0$ and $\kappa >\left\vert \tau
\right\vert $ then $\alpha $ is a part of a spacelike circular helix,

\begin{equation*}
\alpha \left( s\right) =\frac{1}{\left\Vert \omega \right\Vert ^{2}}\left(
\tau \left\Vert \omega \right\Vert s,~\kappa \cos \left[ \left\Vert \omega
\right\Vert s\right] ,~\kappa \sin \left[ \left\Vert \omega \right\Vert s%
\right] \right) .
\end{equation*}

\textbf{Result 3.1.5. }Let $\alpha $ be a spacelike curve with timelike
binormal and timelike $\omega $.$~$If the curve $\alpha $ with $\kappa =$cons%
$\tan $t $>0,~\tau =0$ then $\alpha $ is a part of a circle.\ \ \ \ \ \ \ \
\ \ \ \ \ \ \ \ \ \ \ \ \ \ \ \ \ \ \ \ \ \ \ \ \ \ \ \ \ \ \ \ \ \ \ \ \ \
\ \ \ \ \ \ \ \ \ \ \ \ \ \ \ \ \ \ \ \ \ \ \ \ \ \ \ \ \ \ \ \ \ \ \ \ \ \
\ \ \ \ \ \ \ \ \ \ \ \ \ \ \ \ \ \ \ \ \ \ 

From Lemma 3.1 in \cite{choi}, we can write the following result:

\textbf{Result 3.1.6. }There is no spacelike $W$-curve with timelike
binormal with condition $\left\vert \tau \right\vert =\left\vert \kappa
\right\vert .$\textbf{\ \ \ \ \ \ \ \ \ \ } \ \ \ \ \ \ \ \ \ \ \ \ \ \ \ \
\ \ \ \ \ \ \ \ \ \ \ \ \ \ \ \ \ \ \ \ \ \ \ \ \ \ \ \ \ \ \ \ \ \ \ \ \ \
\ \ \ \ \ \ \ \ \ \ \ \ 

\textbf{Example 3.1.7. }Let $\alpha \left( s\right) =\left( \sinh s,~\cosh
s,~\sqrt{2}s\right) ~$be a unit speed spacelike hyperbolic helix with
timelike binormal and spacelike $\omega $ such that \ \ \ \ \ \ \ \ \ \ \ \
\ \ \ \ \ \ \ \ \ \ \ \ \ \ \ \ \ \ \ \ \ \ \ \ \ \ \ \ \ \ \ \ \ \ \ \ \ \
\ \ \ \ \ \ \ \ \ \ \ \ \ \ \ \ \ \ \ \ \ \ \ \ \ \ \ \ \ \ \ \ \ \ \ \ \ \
\ \ \ \ \ \ \ \ \ \ \ \ \ \ \ \ \ \ \ 

\begin{eqnarray*}
t &=&\left( \cosh s,~\sinh s,~\sqrt{2}\right) \\
n &=&\left( \sinh s,~\cosh s,~0\right) \\
b &=&\left( \sqrt{2}\cosh s,~\sqrt{2}\sinh s,~1\right) ,~\kappa =1~\text{and 
}\tau =\sqrt{2}.
\end{eqnarray*}

If $\alpha $ is a spacelike curve then its involute curve is a spacelike. In
this situation, the involutes of the curve $\alpha $ can be given by the
equation

\begin{equation*}
\alpha ^{\ast }\left( s\right) =\left( \sinh s+\left\vert c-s\right\vert
\cosh s,~\cosh s+\left\vert c-s\right\vert \sinh s,~c\sqrt{2}\right) ,
\end{equation*}%
where $c\in 
\mathbb{R}
.~$One can see a special example of such a curve $\alpha ~$as Fig. 2. and
its involute curve $\alpha ^{\ast }~$as Fig. 3. when $s=\left[ -5,~5\right] $
\ and $c=2.$

\ \ \ \ \ \ \FRAME{itbpFU}{2.5754in}{2.5754in}{0in}{\Qcb{\textbf{Figure 2.}
Spacelike curve $\protect\alpha ^{\ast }$}}{}{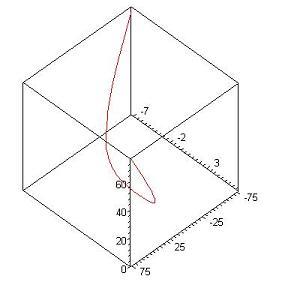}{\special{language
"Scientific Word";type "GRAPHIC";maintain-aspect-ratio TRUE;display
"USEDEF";valid_file "F";width 2.5754in;height 2.5754in;depth
0in;original-width 3in;original-height 3in;cropleft "0";croptop
"1";cropright "1";cropbottom "0";filename 'fig2.jpg';file-properties
"XNPEU";}}\ \FRAME{itbpFU}{2.5858in}{2.5858in}{0in}{\Qcb{\textbf{Figure 3.}
Involute curve $\protect\alpha ^{\ast }$}}{}{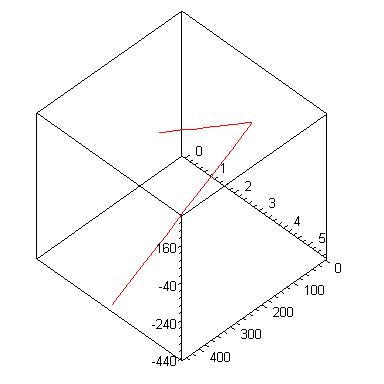}{\special{language
"Scientific Word";type "GRAPHIC";maintain-aspect-ratio TRUE;display
"USEDEF";valid_file "F";width 2.5858in;height 2.5858in;depth
0in;original-width 3in;original-height 3in;cropleft "0";croptop
"1";cropright "1";cropbottom "0";filename 'fig3.jpg';file-properties
"XNPEU";}}\ \ \ \ \ \ \ \ \ \ \ \ \ \ \ \ \ \ \ \ \ \ \ \ \ \ \ \ \ \ \ \ \
\ \ \ \ \ \ \ \ \ \ \ \ \ \ \ \ \ \ \ \ \ \ \ \ \ \ \ \ \ \ \ \ \ \ \ \ \ \
\ \ \ \ \ \ \ \ \ \ \ \ \ \ \ \ \ \ \ ~\ \ \ \ \ \ \ \ \ \ \ \ \ \ \ \ \ \ \
\ \ \ \ \ \ \ \ \ \ \ \ \ \ \ \ \ \ \ \ \ \ \ \ \ \ \ \ \ \ \ \ \ \ \ \ \ \
\ \ \ \ \ \ \ \ \ \ \ \ \ \ \ \ \ \ \ \ \ \ \ \ \ \ \ \ \ \ \qquad \qquad
\qquad\ \ \ \ \ \ \ \ \ \ \ \ \ \ \ \ \ \ \ \ \ \ \ \ \ \ \ \ \ \ \ \ \ \ \
\ \ \ \ \ \ \ \ \ \ \ \ \ \ \ \ \ \ \ \ \ \ \ \ \ \ \ \ \ \ \ \ \ \ \ \ \ \
\ \ \ \ \ \ \ \ \ \ \ \ \ \ \ \ \ \ \ \ \ \ \ \ \ 

\bigskip The short calculations give the following equation of the spherical
indicatrices of the involute curve $\alpha ^{\ast }$ and its natural lifts.
\ \ \ \ \ \ \ \ \ \ 

\begin{center}
$%
\begin{array}{c}
\text{ \ \ \ \ \ \ \ \ \ \ \ \ }\alpha _{t^{\ast }}^{\ast }=t^{\ast }=\left(
\sinh s,~\cosh s,~0\right) \text{\ } \\ 
\text{ \ \ \ \ \ \ \ \ \ \ \ \ }\alpha _{n^{\ast }}^{\ast }=n^{\ast }=\left(
\cosh s,~\sinh s,~0\right) \\ 
\alpha _{b^{\ast }}^{\ast }=b^{\ast }=\left( 0,~0,~1\right)%
\end{array}%
\begin{array}{c}
\text{ \ \ \ \ \ \ \ \ \ \ \ }\overline{\alpha ^{\ast }}_{t^{\ast }}=\left(
\cosh s,~\sinh s,~0\right) \\ 
\text{\ \ \ \ \ \ \ \ \ \ \ \ \ }\overline{\alpha ^{\ast }}_{n^{\ast
}}=\left( \sinh s,~\cosh s,~0\right) \\ 
\overline{\alpha ^{\ast }}_{b^{\ast }}=\left( 0,~0,~0\right)%
\end{array}%
$
\end{center}

Since 
\begin{equation*}
g\left( \alpha _{t^{\ast }}^{\ast \prime },~\alpha _{t^{\ast }}^{\ast \prime
}\right) =1>0
\end{equation*}%
$\alpha _{t^{\ast }}^{\ast }$ is spacelike. For being $\alpha _{t^{\ast
}}^{\ast }$ is a spacelike curve, its spherical image is geodesic which lies
on the Lorentzian unit sphere $S_{1}^{2}~$as Fig. 4. and natural lift curve
of the tangent indicatrix as Fig. 5. One consider the principal normal
indicatrix is a geodesic which lies on $H_{0}^{2}$ as Fig. 6 and its natural
lift as Fig. 7.

\ \ \ \ \ \ \ \FRAME{dtbpFU}{2.4811in}{2.7354in}{0pt}{\Qcb{\textbf{Figure 4.}
Spherical image of tangent indicatrix of the involute curve $\protect\alpha %
^{\ast }$}}{}{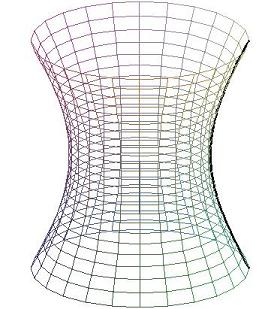}{\special{language "Scientific Word";type
"GRAPHIC";maintain-aspect-ratio TRUE;display "USEDEF";valid_file "F";width
2.4811in;height 2.7354in;depth 0pt;original-width 2.917in;original-height
3.2188in;cropleft "0";croptop "1";cropright "1";cropbottom "0";filename
'fig4.jpg';file-properties "XNPEU";}}\ \ \ \FRAME{dtbpFU}{2.2822in}{2.7363in%
}{0pt}{\Qcb{\textbf{Figure 5.} Tangent indicatrix of the involute curve $%
\protect\alpha ^{\ast }$\textbf{\ }and its natural lift}}{}{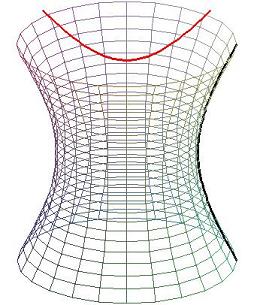}{%
\special{language "Scientific Word";type "GRAPHIC";maintain-aspect-ratio
TRUE;display "USEDEF";valid_file "F";width 2.2822in;height 2.7363in;depth
0pt;original-width 2.6455in;original-height 3.1773in;cropleft "0";croptop
"1";cropright "1";cropbottom "0";filename 'fig5.jpg';file-properties
"XNPEU";}}

\ \ \ \ \ \ \FRAME{dtbpFU}{2.5962in}{2.5962in}{0pt}{\Qcb{\textbf{Figure 6.}
Spherical image of principal normal indicatrix of the involute curve $%
\protect\alpha ^{\ast }$}}{}{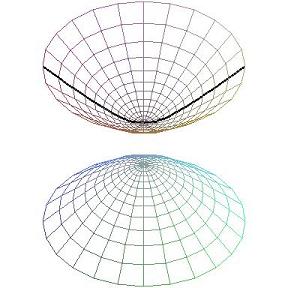}{\special{language "Scientific
Word";type "GRAPHIC";maintain-aspect-ratio TRUE;display "USEDEF";valid_file
"F";width 2.5962in;height 2.5962in;depth 0pt;original-width
3in;original-height 3in;cropleft "0";croptop "1";cropright "1";cropbottom
"0";filename 'fig6.jpg';file-properties "XNPEU";}} \ \ \ \FRAME{dtbpFU}{%
2.6066in}{2.6066in}{0pt}{\Qcb{\textbf{Figure 7.} Principal normal indicatrix
of the involute curve $\protect\alpha ^{\ast }$\textbf{\ }and its natural
lift}}{}{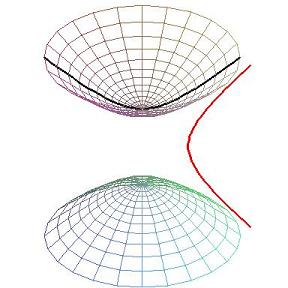}{\special{language "Scientific Word";type
"GRAPHIC";maintain-aspect-ratio TRUE;display "USEDEF";valid_file "F";width
2.6066in;height 2.6066in;depth 0pt;original-width 3in;original-height
3in;cropleft "0";croptop "1";cropright "1";cropbottom "0";filename
'fig7.jpg';file-properties "XNPEU";}}

\subsection{\textbf{The natural lift} \textbf{of principal normal indicatrix}
\textbf{of the} \textbf{curve }$\protect\alpha ^{\ast }$}

Let $\alpha ~$be a spacelike curve with timelike binormal and spacelike $%
\omega $ $\left( \left\vert \kappa \right\vert <\left\vert \tau \right\vert
\right) $. In this section, we will investigate how $\alpha $ must be a
curve satisfying the condition that the natural lift curve $\overline{\alpha
^{\ast }}_{n^{\ast }}~$of \ $\alpha _{n^{\ast }}^{\ast }$ is an integral
curve of geodesic spray, where $\alpha _{n^{\ast }}^{\ast }$ is the
spherical indicatrix of principal normal vector of $\alpha ^{\ast }$. If the
natural lift curve $\overline{\alpha ^{\ast }}_{n^{\ast }}$\ is an integral
curve of the geodesic spray, then by means of Lemma 2.1. we have

\begin{equation}
\overline{\overline{D}}_{\alpha _{n^{\ast }}^{\prime \ast }}\alpha _{n^{\ast
}}^{\prime \ast }=0,  \label{u5}
\end{equation}%
and from the Lemma 2.2. I. and the equation (5) we get,

\begin{equation*}
\left[ \left( \sigma ^{\prime }\cosh \theta +\theta ^{\prime }\sigma \sinh
\theta -\frac{\kappa }{k_{n}}\right) t+\left( -\frac{k_{n}^{\prime }}{%
k_{n}^{2}}\right) n+\left( \frac{\tau }{k_{n}}-\sigma ^{\prime }\sinh \theta
-\theta ^{\prime }\sigma \cosh \theta \right) b\right] \frac{1}{\left\Vert
\omega \right\Vert k_{n}}=0,
\end{equation*}%
where $\sigma =\frac{\gamma _{n}}{k_{n}}~$( $\gamma _{n}=\frac{\theta
^{\prime }}{\left\Vert \omega \right\Vert }~$and $k_{n}=\frac{1}{\left\Vert
\omega \right\Vert }\sqrt{\theta ^{\prime ^{2}}+\left\Vert \omega
\right\Vert ^{2}}$ are the geodesic curvatures of the curve $\alpha $ with
respect to $S_{1}^{2}$ and $\hbox{\bf E}_{1}^{3}$, respectively.) and $%
\overline{\overline{D}}~$is the connection of hyperbolic sphere $H_{0}^{2}$.
Since $\left\{ t,~n,~b\right\} $ are linear independent, we get

\begin{eqnarray*}
\sigma ^{\prime }\cosh \theta +\theta ^{\prime }\sigma \sinh \theta -\frac{%
\kappa }{k_{n}} &=&0 \\
\frac{k_{n}^{\prime }}{k_{n}^{2}} &=&0 \\
\frac{\tau }{k_{n}}-\sigma ^{\prime }\sinh \theta -\theta ^{\prime }\sigma
\cosh \theta &=&0,
\end{eqnarray*}%
and we obtain

\begin{equation*}
\gamma _{n}=\text{cons}\tan \text{t, }k_{n}=\text{cons}\tan \text{t. }
\end{equation*}%
Therefore, we can write the following result.

\textbf{Result 3.2.1. }If the geodesic curvatures of the evolute curve $%
\alpha $ with respect to $S_{1}^{2}$ and $\hbox{\bf E}_{1}^{3}$ are
constant, then the spherical indicatrix $\alpha _{n^{\ast }}^{\ast }$ is a
geodesic on the hyperbolic sphere $H_{0}^{2}$. In this case, the natural
lift $\overline{\alpha ^{\ast }}_{n^{\ast }}$ of $\alpha _{n^{\ast }}^{\ast
} $ is an integral curve of the geodesic spray on the tangent bundle $%
T\left( H_{0}^{2}\right) $. In particular, if the evolute curve $\alpha $ is
a spacelike curve with timelike binormal and timelike $\omega ~\left(
\left\vert \kappa \right\vert >\left\vert \tau \right\vert \right) $ , then
the similar result can be easily obtained by taking $S_{1}^{2}~$instead of $%
H_{0}^{2}$ in following same procedure. \ \ \ \ \ \ \ \ \ \ \ \ \ \ \ \ \ \
\ \ \ \ \ \ \ \ \ \ \ \ \ \ \ \ \ \ \ \ \ \ \ \ \ \ \ \ \ \ \ \ \ \ \ \ \ \
\ \ \ \ \ \ \ \ \ \ \ \ \ \ \ \ \ \ \ \ \ \ \ \ \ \ \ \ \ \ \ \ \ \ \ \ \ \
\ \ \ \ \ \ \ \ \ \ \ \ \ \ \ \ \ \ \ \ \ \ \ \ \ \ \ \ \ \ \ \ \ \ \ \ \ \
\ \ \ \ \ \ \ \ \ \ 

\subsection{\textbf{The natural lift} \textbf{of binormal indicatrix} 
\textbf{of the} \textbf{curve }$\protect\alpha ^{\ast }$}

Let $\alpha ~$be a spacelike curve with timelike binormal and spacelike $%
\omega $ $\left( \left\vert \kappa \right\vert <\left\vert \tau \right\vert
\right) $. We will investigate how $\alpha $ must be a curve satisfying the
condition that the natural lift curve $\overline{\alpha ^{\ast }}_{b^{\ast
}}~$is an integral curve of geodesic spray, where $\alpha _{b^{\ast }}^{\ast
}$ is the spherical indicatrix of binormal vector of $\alpha ^{\ast }~$and $%
\overline{\alpha ^{\ast }}_{b^{\ast }}~$is the natural lift of the curve $%
\alpha _{b^{\ast }}^{\ast }$. If the natural lift curve $\overline{\alpha
^{\ast }}_{b^{\ast }}$\ is an integral curve of the geodesic spray, then by
means of Lemma 2.1. we have

\begin{equation}
\overline{D}_{\alpha _{b^{\ast }}^{\prime \ast }}\alpha _{b^{\ast }}^{\prime
\ast }=0,  \label{u6}
\end{equation}%
from the Lemma 2.2. I. and the equation (6) we have,

\begin{equation*}
\frac{\left\Vert \omega \right\Vert }{\theta ^{\prime }}n=0.
\end{equation*}%
Since $\left\{ t,~n,~b\right\} $ are linear independent, we obtain

\begin{equation*}
\kappa =0,~\tau =0.
\end{equation*}%
Thus, we can give the following result.

\textbf{Result 3.3.1. }The spherical indicatrix $\alpha _{b^{\ast }}^{\ast }$
of the involute curve $\alpha ^{\ast }$ can not be a geodesic line on the
Lorentzian sphere $S_{1}^{2}$, because, the evolute curve $\alpha $ whose
curvature and torsion are equal to $0$ is a straight line. In this case $%
\left( \alpha ,~\alpha ^{\ast }\right) $ can not occur the involute-evolute
curve couple. Therefore, the natural lift $\overline{\alpha ^{\ast }}%
_{b^{\ast }}$ of the curve $\alpha _{b^{\ast }}^{\ast }$ can never be an
integral curve of the geodesic spray on the tangent bundle $T\left(
S_{1}^{2}\right) $. If the evolute curve $\alpha $ is a spacelike curve with
timelike binormal and timelike $\omega $ , then the similar result can be
easily obtained by taking $S_{1}^{2}~$instead of $H_{0}^{2}$ in following
same procedure. \ \ \ \ \ \ \ \ \ \ \ \ \ \ \ \ \ \ \ \ \ \ \ \ \ \ \ \ \ \
\ \ \ \ \ \ \ \ \ \ \ \ \ \ \ \ \ \ \ \ \ \ \ \ \ \ \ \ \ \ \ \ \ \ \ \ \ \
\ \ \ \ \ \ \ \ \ \ \ \ \ \ \ \ \ \ \ \ \ \ \ \ \ \ \ \ \ \ \ \ \ \ \ \ \ \
\ \ \ \ \ \ \ \ \ \ \ \ \ \ \ \ \ \ \ \ \ \ \ \ \ \ \ \ \ \ \ \ \ \ \ \ \ \
\ \ \ \ \ \ \ \ \ \ \ \ \ \ \ \ \qquad \qquad \qquad \qquad


\begin{thebibliography}{99}
\bibitem{ali1} A.T. Ali, Position vectors of spacelike general helices in
Minkowski 3-space, Nonlinear Anal. TMA, 73 (2010) 1118-1126.

\bibitem{ali3} A.T. Ali, R. Lopez, Slant helices in Minkowski space , J.
Korean Math. Soc. 48 (2011) 159-167.

\bibitem{barr} M. Barros, A. Ferrandez, P. Lucas, M.A. Merono, General
helices in the three-dimensional Lorentzian space forms, Rocky Mountain J.
Math. 31 (2001) 373-388.

\bibitem{bili} M. Bilici, M. \c{C}al\i \c{s}kan, On the involutes of the
spacelike curve with a timelike binormal in Minkowski 3-space, Int. Math.
Forum. 4 (2009) 1497-1509.

\bibitem{bilici1} M. Bilici, M. \c{C}al\i \c{s}kan, \.{I}. Aydemir, The
Natural Lift Curves and the Geodesic Sprays for the Spherical Indicatrices
of the Pair of Evolute-Involute Curves, Int. J. of Appl. Math. 11 (2003)
415-420.

\bibitem{bilici2} M. Bilici, Natural lift curves and the geodesic sprays for
the spherical indicatrices of the involutes of a timelike curve in Minkowski
3-space, Int. J. Phys Sci. 6 (2011) 4706-4711.

\bibitem{choi} J.H. Choi, Y.H. Kim, A.T. Ali, Some associated curves of
Frenet non-lightlike curves in $\hbox{\bf E}_{1}^{3}$, J. Math. Anal. Appl.
394 (2012) 712-723.

\bibitem{cal} M. \c{C}al\i \c{s}kan, M. Bilici, Some Characterizations for
the Pair of Involute-Evolute Curves in Euclidean Space \textit{$\hbox{\bf E}%
^{3}$}, Bull. of Pure and Appl.Sci. 21 (2002) 289-294.

\bibitem{call} M. \c{C}al\i \c{s}kan, A.\.{I}. Sivrida\u{g}, H.H. Hac\i
saliho\u{g}lu, Some Characterizations For The Natural Lift Curves and The
Geodesic Spray, Commun. Fac. Sci. Univ. 33 (1984) 235-242.

\bibitem{Do Carmo} Do Carmo MP, Differantial Geometry of Curves and
Surfaces, Prentice Hall, Englewood Cliffs, NJ, 1976.

\bibitem{hs} Hac\i saliho\u{g}lu HH, Differantial Geometry, Ankara
University Faculty of Science Press, Ankara, 2000.

\bibitem{ilarslan} K. \.{I}larslan, \"{O}. Boyac\i o\u{g}lu, Position
vectors of a timelike and a null helix in Minkowski 3-space, Chaos, Solitons
and Fractals, 38 (2008) 1383--1389

\bibitem{mill} R.S. Millman, G.D. Parker, Elements of Differential Geometry.
Prentice-Hall Inc., Englewood Cliffs, New Jersey, 1977.

\bibitem{Onei} B. O'Neill, Semi-Riemannian Geometry with Application to
relativity, Academic Press, New York, 1983.

\bibitem{rat} J.G. Ratcliffe, Foundations of Hyperbolic Manifolds,
Springer-Verlag New York, Inc., New York, 1994.

\bibitem{thorpe} J.A. Thorpe, Elemantary Topics In Differantial Gemetry,
Springer-Verlag, New York, Heidelberg-Berlin, 1979.

\bibitem{ugur} H.H. U\u{g}urlu, On The Geometry of Time-like Surfaces,
Commun. Fac. Sci. Univ. Ank. Series A1, 46 (1997) 211-223.

\bibitem{walrave} J. Walrave, Curves and Surfaces in Minkowski Space. PhD
thesis, K.U. Leuven, Faculty of Science, Leuven, 1995.

\bibitem{Kula} L. Kula, Y. Yayl\i , On slant helix and its spherical
indicatrix. Appl. Math. and Comput. 169 (2005) 600-607.

\bibitem{Suha} S. Y\i lmaz, E. \"{O}zy\i lmaz, Y. Yayl\i , M. Turgut,
Tangent and trinormal spherical images of a time-like curve on the
pseudohyperbolic space. Proc. Est. Acad. Sci. 59 (2010.) 216--224.

\bibitem{Esen} E. \.{I}yig\"{u}n, The tangent spherical image and ccr-curve
of a time-like curve in \textit{$\hbox{\bf L}^{3}$}, J. Inequal. Appl.
(2013) doi:10.1186/1029-242X-2013-55.
\end{thebibliography}
\end{document}